\documentstyle[12pt,amssymb,epsf]{article}

\def\R {{\Bbb R}}
\def\N {{\Bbb N}}

\def\t {\theta}
\def\a {\alpha}

\newtheorem{theorem}{Theorem}
\newtheorem{proposition}{Proposition}
\newtheorem{definition}{Definition}
\newtheorem{lemma}{Lemma}
\newtheorem{remark}{Remark}

\begin{document}
\title{Limit of Karcher's Saddle towers}
\author{M. Magdalena Rodr\'\i guez}

\maketitle

\section{Introduction}\label{secintro}
In 1835, Scherk~\cite{sche1} showed a singly periodic minimal surface
$S$ in $\R^3$, which may be viewed as the desingularization of two
vertical planes meeting at a right angle.  This surface $S$ was
generalized later on to a one-parameter family of singly periodic
minimal surfaces $S_\t$ in $\R^3$, where $\t\in(0,\frac{\pi }{2}]$ is
the angle between the asymptotic vertical planes (in particular,
$S= S_\frac{\pi }{2}$).  In the quotient by its shortest period vector,
each $S_\t$ has genus zero and four ends asymptotic to flat vertical
annuli. Annular ends of this kind are called {\it Scherk-type ends}.
These singly periodic Scherk minimal surfaces have recently been
classified in~\cite{mrw1} as the only properly embedded singly
periodic minimal surfaces with four Scherk-type ends in the quotient.

H. Karcher~\cite{ka4} generalized the previous Scherk minimal surfaces by
constructing, for each natural $n\geq 2$, a $(2n-3)$-parameter family
of singly periodic minimal surfaces with genus zero and $2n$
Scherk-type ends in the quotient. These surfaces, called {\it saddle
  towers}, are the only properly embedded singly periodic minimal
surfaces in $\R^3$ with genus zero and finitely many Scherk-type ends
in the quotient, see~\cite{PeTra1}. Note that for $n=2$ we obtain the
singly periodic Scherk minimal surfaces.

Let us now recall the construction of the saddle towers: consider any
convex polygonal domain $\Omega_n$ whose boundary consists of $2n$
edges of length one, with $n\geq 2$, and mark its edges alternately by
$\pm\infty$.  Assume $\Omega_n$ is non-special (see
definition~\ref{def1} below).  By a theorem of Jenkins and
Serrin~\cite{jes1}, there exists a function $u_n$ which solves the
Jenkins-Serrin problem on $\Omega_n$; i.e.  $u_n$ is a minimal graph
defined on $\Omega_n$ which diverges to $\pm\infty$, as indicated by
the marking, when we approach to the edges of $\Omega_n$.  The
boundary of this minimal graph consists of $2n$ vertical lines above
the vertices of $\Omega_n$. Hence the conjugate minimal surface of
this graph is bounded by $2n$ horizontal symmetry curves, lying in two
horizontal planes at distance 1 from each other.  By reflecting about
one of the two symmetry planes, we obtain a fundamental domain for a
saddle tower $M_n$ with period $T=(0,0,2)$ and $2n$ Scherk-type ends
in the quotient.

\begin{definition}
\label{def1}
We say that a convex polygonal domain with $2n$ unitary edges is special if
$n\geq 3$ and its boundary is a parallelogram with two sides of length
one and two sides of length $n-1$.
\end{definition}

\begin{remark}[\cite{mrt}]
The bounded convex polygonal domains with edges of length one which
fail to satisfy the hypothesis of the theorem of Jenkins and
Serrin~\cite{jes1} are precisely the special domains.
\end{remark}

In this paper we study the possible limits of saddle towers by taking
limits of sequences of minimal graphs $u_n$ as above. We normalize so
that the segment of vertices $(0,0),(1,0)$ is one of the edges of the
convex polygonal domain $\Omega_n$
where $u_n$ is defined, and such edge is marked by $+\infty$.

\begin{theorem}
\label{th}
Let $M_n\subset\R^3$ be a saddle tower with $2n$ Scherk-type ends in
the quotient, and $\Omega_n$ be its associated normalized convex
polygonal domain in the above construction.  Suppose $\{\Omega_n\}_n$
does not converge to a straight line nor half a straight line.  Then a
subsequence of $\{M_n\}_n$ converges uniformly on compact sets of
$\R^3$ with multiplicity one to either one of the saddle towers with
infinitely many ends described in~\cite{mrt}, any singly periodic
Scherk minimal surface, a doubly periodic Scherk minimal surface of
angle $\frac{\pi}{2}$, or a KMR example $M_{\t,\a,0}$ studied
in~\cite{mrod1} (also called {\it toroidal halfplane layer} by
Karcher~\cite{ka4}).
\end{theorem}

\section{Preliminaries}\label{secprelim}
In this section we present some general results for minimal graphs
explained in~\cite{mrt} and based on the ideas of Jenkins and Serrin,
developed by Collin and Mazet.

Let $u=u(x_1,x_2)$ be a solution of the minimal graph equation,
\begin{equation}\label{eqmin}
\mbox{div}\left(\frac{\nabla u}{\sqrt{1+|\nabla u|^2}}\right)=0 ,
\end{equation}
defined on a domain $\Omega\subset\R^2$.  By an elementary computation,
we obtain that the form $d\psi_u=\frac{u_{x_1}}{\sqrt{1+|\nabla
    u|^2}}\, d x_2-\frac{u_{x_2}}{\sqrt{1+|\nabla u|^2}}\, d x_1$
is closed. Hence it defines a function $\psi_u=\psi_u(x_1,x_2)$,
called {\it conjugate function of $u$}, which is well defined up to an
additive constant. 
In fact, $\psi_u$ coincides with the third  coordinate function of the
conjugate minimal surface of the graph of $u$, written as a function on
the $(x_1,x_2)$-parameters (although such conjugate surface does not
coincide with the graph of $\psi_u$).
It is straightforward to check that $\psi_u$ is a Lipschitz function,
in particular it can be extended continuously to $\partial\Omega$.
Moreover, the following lemma holds.

\begin{lemma}
  \label{lem1}
  Given a solution $u$ of (\ref{eqmin}) on a domain $\Omega\subset\R^2$, we have:
  \begin{enumerate}
  \item[(i)] For every domain $D\subset\Omega$,
    $\int_{\partial D} d\psi_u=0$.
  \item[(ii)] Let $T\subset\partial\Omega$ be a bounded arc oriented as
    $\partial D$.  Then,
    \[
    \left|\int_T d\psi_u\right|\leq |T|,
    \]
    and $\int_T d\psi_u= |T|$ (resp. $-|T|$) if and only if $u$
    diverges to $+\infty$ (resp. to $-\infty$) as one approaches $T$
    within $\Omega$, in which case $T$ must be a straight segment.
  \end{enumerate}
\end{lemma}

Let $u_n$ be a solution of~(\ref{eqmin}) on a domain $\Omega_n\subset\R^2$. We
define the {\it limit domain} $\Omega_\infty$ of the domains $\Omega_n$ as the set of
points in $\R^2$ that admit a neighborhood contained in every
$\Omega_n$, for $n$ large enough, and say that $\{\Omega_n\}_n$
converges to $\Omega_\infty$. 
Consider the {\it convergence domain} of $\{u_n\}_n$, defined as
\[
{\cal B}(u_n)=\left\{p\in\Omega_\infty\ |\ \{|\nabla u_n|(p)\}_n
  \mbox{ is bounded }\right\}.
\]
For each component $D$ of ${\cal B}(u_n)$ and any point $q\in D$, there is a subsequence of
$\{u_n-u_n(q_n)\}_n$ converging uniformly on compact sets of $D$ to a
solution of~(\ref{eqmin}), where $q_n\in\Omega_n$ with $q_n\to q$.
Moreover, 
\[
\Omega_\infty-{\cal B}(u_n)=\cup_{i\in I} L_i ,
\]
where each $L_i$, called a {\it divergence line}, is a component of
the intersection of a straight line with~$\Omega_\infty$.  Clearly, to
ensure the convergence of a subsequence of the vertical translated
$u_n$ on $\Omega_\infty$, it suffices to prove there are no divergence
lines.

\begin{lemma}\label{lem2}
  Let $\{u_n\}_n$ be a sequence of minimal graphs as above, and denote
  by $\psi_n$ the conjugate function of $u_n$, for every $n\in\N$.
\begin{enumerate}
\item[(i)] Let $T$ be a straight segment contained in a divergence
  line, $T\subset\Omega_n$ for $n$ big enough.
  Then $\int_T d\psi_n\to\pm|T|$.
\item[(ii)] A divergence line cannot finish 
  at an interior point of an open straight segment
  $T\subset\partial\Omega_\infty$, if for each $n$ there exists a
  straight segment $T_n$ in $\Omega_n$ such that $u_n$ diverges to
  $+\infty$ when we approach $T_n$ within $\Omega_n$ and $T_n\to T$.
\end{enumerate}
\end{lemma}

Finally, we have the following uniqueness result for the limit $u$
under some constraints.
\begin{lemma}[\cite{mazet3}]
\label{lem5}
Let $u,v$ be two solutions of (\ref{eqmin}) in a
domain~$\Omega$, whose conjugate functions $\psi_u,\psi_v$ are
bounded in $\Omega$ and coincide in $\partial\Omega$.  Then
$u-v$ is constant in $\Omega$.
\end{lemma}

\section{Taking limits of saddle towers}
For every $n\in\N$, let $\Omega_n\subset\R^2$ be a non-special,
convex, bounded polygonal domain with $2n$ unitary edges.  Denote its
vertices as $p_i^n$, $i=0,\cdots,2n-1$, by following the natural
cyclic ordering induced by the positive orientation of
$\partial\Omega_n$.  Assume that $p_0^n=(0,0)$ and $p_1^n=(1,0)$ for
every $n\in\N$.  After passing to a subsequence, the sequence of
domains $\{\Omega_n\}_n$ converges to either a straight line, half a
straight line, or a convex, unbounded, polygonal domain
$\Omega_\infty$ with unitary edges, whose vertices are obtained as
limits of the vertices of the domains $\Omega_n$.  Assume that the two
first cases do not happen.

Let $u_n$ be the solution to the minimal graph equation~(\ref{eqmin})
on $\Omega_n$ which takes boundary values $+\infty$ on the edges
$(p_{2i}^n,p_{2i+1}^n)$, and $-\infty$ on the edges
$(p_{2i-1}^n,p_{2i}^n)$, $i=0,\cdots, n-1$.  Our aim consists of
studying the possible limits for $\{u_n\}_n$.

We denote by $\psi_n$ the conjugate function of $u_n$ verifying
$\psi_n(p_0^n)=0$.  From Lemma~\ref{lem1}-{\it (ii)} we have
$\int_{p_i^n}^{p_{i+1}^n} d\psi_n=(-1)^i$, which implies that
$\psi_n(p_i^n)$ is equal to $0$ if $i$ is even, and equal to $1$ when
$i$ is odd.  Moreover, $\psi_n$ is an affine function on each edge of
$\Omega_n$, so $0\leq\psi_n\leq 1$ on $\partial\Omega_n$.  And by the
maximum principle, $0\leq\psi_n\leq 1$ in $\Omega_n$.

Finally, we will say that a vertex of $\Omega_\infty$ is {\it even}
(resp. {\it odd}) if it is obtained as limit of vertices $p_{2i}^n$
(resp. $p_{2i+1}^n$) of the domains $\Omega_n$, where $i=0,\cdots,
2n-1$ (this is possible because we have fixed the vertices
$p_0^n,p_1^n$).  In particular, $p_0=(0,0)$ is an even vertex of
$\Omega_\infty$, and $p_1=(1,0)$ is an odd vertex.

\begin{proposition}
  \label{prop}
  If the distance between any two vertices of $\Omega_\infty$ with
  different parity is strictly bigger than one, and $q_n\in\Omega_n$
  with $q_n\to q\in\Omega_\infty$, then $\{u_n-u_n(q_n)\}_n$ converges to a
  minimal graph $u_\infty$ with boundary values $\pm\infty$ disposed
  alternately on $\partial\Omega_\infty$ and whose conjugate graph
  lies in the horizontal slab $\{(x_1,x_2,x_3)\ |\ 0\leq x_3\leq 1\}$.
\end{proposition}
\begin{proof}
  Let us prove there are no divergence lines for $\{u_n\}_n$ in the
  setting of Proposition~\ref{prop}. Suppose by contradiction that
  there exists one such divergence line $L$. We deduce from
  Lemma~\ref{lem2}-{\it (i)} that $L$ must have length no bigger than
  one, since $0\leq\psi_n\leq 1$ for every $n\in\N$. Thus
  Lemma~\ref{lem2}-{\it (ii)} says that $L$ must be a segment (of
  length at most one) joining two different vertices $p_i,p_j$ of
  $\Omega_\infty$.  Let $p_i^n,p_j^n$ be vertices of $\Omega_n$ such
  that $p_i^n\to p_i$ and $p_j^n\to p_j$.  Since
  $|\int_{p_i^n}^{p_j^n}d\psi_n|\to|p_i-p_j|$ by using
  Lemma~\ref{lem2}-{\it (i)}, but
  $|\int_{p_i^n}^{p_j^n}d\psi_n|=|\psi_n(p_j^n)-\psi_n(p_i^n)|$ can
  only equal $0$ or $1$, we deduce that the only possibility is
  $|p_i-p_j|=1$, and so $|\psi_n(p_j^n)-\psi_n(p_i^n)|=1$ for $n$
  large. In particular, the vertices $p_i,p_j$ have different parity
  and satisfy $|p_i-p_j|=1$, which contradicts the hypothesis in
  Proposition~\ref{prop}.
  
  Hence there exists a subsequence of $\{u_n-u_n(q_n)\}_n$ converging on
  compact subsets of $\Omega_\infty$ to a minimal graph $u_\infty$.
  Moreover, we deduce from Lemma~\ref{lem1}-{\it (ii)} that $u_\infty$
  takes boundary values $\pm\infty$ alternately on the unitary edges in
  $\partial\Omega_\infty$.  Since the conjugate function $\psi_\infty$
  of $u_\infty$ can be obtained as limit of $\{\psi_n\}_n$, we deduce that
  $0\leq \psi_\infty\leq 1$ on $\Omega_\infty$.

  Finally, we obtain from Lemma~\ref{lem5} that it is not only a
  subsequence but the whole sequence $\{u_n-u_n(q_n)\}_n$ what
  converges to $u_\infty$, and Proposition~\ref{prop} is proven.
\end{proof}

\begin{remark}
  \label{rem}
  From Lemma~\ref{lem5}, we know that the graph $u_\infty$ obtained in
  Proposition~\ref{prop} must be half a singly periodic Scherk minimal
  surface of angle $\frac{\pi}{2}$, if $\Omega_\infty$ is a halfplane; one of the graphs
  studied in~\cite{mrod2} (i.e. a piece of a KMR example
  $M_{\t,\a,\frac{\pi}{2}}$), if the limit domain $\Omega_\infty$ is a
  strip; or one of the graphs constructed in~\cite{mrt}, in other
  case.
\end{remark}

Now suppose we are not in the setting of Proposition~\ref{prop}; this
is, suppose there exist two vertices $p_i,p_j$ of $\Omega_\infty$ with different parity 
such that $0\neq|p_i-p_j|\leq 1$.
It is not very difficult to check that $\Omega_\infty$ must be either a
strip or a special, unbounded, convex polygonal domain, see
Definition~\ref{def2} below.

\begin{figure}\begin{center}
\epsfysize=2cm \epsffile{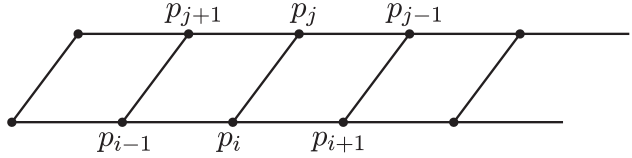}
\end{center}
\caption{An example of special, unbounded, convex polygonal domain.}
\label{fig}
\end{figure}

\begin{definition}
  \label{def2}
  An unbounded convex polygonal domain is said to be special when its
  boundary is made of two parallel half lines and one edge of length
  one (such a domain may be seen as a limit of special domains with
  $2n$ edges, when $n\to\infty$), see Figure~\ref{fig}.
\end{definition}

Let $p_{i_n}^n,p_{j_n}^n$ be two vertices of $\Omega_n$ such that
$p_{i_n}^n\to p_i$ and $p_{j_n}^n\to p_j$. Assume $p_i$ is an even
vertex and $p_j$ is an odd vertex. Then $\psi_n(p_{i_n}^n)=0$ and
$\psi_n(p_{j_n}^n)=1$, so $|\int_{p_{i_n}^n}^{p_{j_n}^n}d\psi_n|=1$.
Thus we deduce from Lemma~\ref{lem1}-{\it (ii)} that the straight
segment $(p_i,p_j)$ is a divergence line.  Moreover, we can similarly
prove that, if $p_{i_n+1}^n, p_{j_n-1}^n$ (resp.  $p_{i_n-1}^n,
p_{j_n+1}^n$) are not consecutive vertices and for each index~$\a$ we
denote by $p_\a$ the vertex of $\Omega_\infty$ so that $p_{\a_n}^n\to
p_\a$, then the straight segment $(p_{i+1},p_{j-1})$ (resp. $(p_{i-1},
p_{j+1})$) is a divergence line.  Following this reasoning we obtain
that the convergence domain ${\mathcal B}(u_n)$ consists of
consecutive translated rhombi with unitary edges (see
Figure~\ref{fig}). By uniqueness we know that, after a suitable
vertical translation, the graphs $u_n$ converge on each such rhombus
to a fundamental piece of a doubly periodic Scherk minimal surface. We
then obtain the following lemma.\\

\begin{lemma}
  \label{lem}
  If $\Omega_\infty$ is a special domain, the convergence domain of
  $\{u_n\}_n$ consists of consecutive translated rhombi
  $R_1,R_2,\dots$ with unitary edges. Moreover, given $q\in R_k$ for
  some $k=1,2,\dots$, the sequence $\{u_n-u_n(q)\}_n$ converges
  uniformly on compact sets of $R_k$ to a fundamental piece of a doubly periodic Scherk minimal
  surface.
\end{lemma}
  
For each $n\in\N$, let $M_n$ be the saddle tower obtained from the
graph~$u_n$. Translate each $M_n$ so that it contains the origin of
$\R^3$, and let us now prove that $\{M_n\}_n$ converges with finite
multiplicity on compact sets of $\R^3$ to a minimal surface
$M_\infty$.  Since we are assuming that $\Omega_\infty$ is not half a
straight line nor a straight line, there exists a uniform radius
$r_0>0$ such that, for each vertex $p_i^n$ of $\Omega_n$, the disk
$D(p_i^n,r_0)$ intersects $\partial\Omega_n$ only along its two edges
with common endpoint $p_i^n$.  We can then prove as in~\cite{mrt},
Section~4, that there exists a constant~$C$ (independent of $n$) such
that the Gauss curvature of $M_n$ is bounded by~$C$.  By the Regular
Neighborhood Theorem, or Rolling Lemma~\cite{ros5,mr7}, $M_n$ has an
embedded tubular neighborhood of radius $1/\sqrt{C}$.  In particular,
we have local area bounds (more precisely, the area of $M_n$ inside
balls of radius $1/\sqrt{C}$ is bounded by some constant).  By
standard result, a subsequence of $\{M_n\}_n$ converges with finite
multiplicity on compact subsets of $\R^3$ to a minimal surface
$M_\infty$.

From Proposition~\ref{prop}, Remark~\ref{rem} and Lemma~\ref{lem}, we
obtain that $M_\infty$ must be the conjugate surface of the
corresponding limit of the graphs $u_n$, i.e.  a KMR example
$M_{\t,\a,0}$, a singly or doubly periodic Scherk minimal surface, or
one of the examples constructed in~\cite{mrt}, which are singly
periodic minimal surfaces with genus zero and one limit end in the
quotient by all their periods. Furthermore, sin $M_\infty$ is not a
plane, the multiplicity of convergence is one.  This completes the
proof of Theorem~\ref{th}.

\begin{remark}
  When the domains $\Omega_n$ converges to either a straight line or
  half a straight line, it may be proven that the Gauss curvature of
  the saddle towers $M_n$ blows-up. And, after scaling to have
  curvature estimates, the graphs $u_n$ converge to half a helicoid
  (see~\cite{mrod2} for a description of this limit).
  In fact, this is the same helicoidal limit we obtain by taking
  limits from fundamental pieces of doubly periodic Scherk minimal
  surfaces. 
  Hence the scaled saddle towers converge to a catenoid.
\end{remark}

\bibliographystyle{plain}

\begin{thebibliography}{10}

\bibitem{jes1}
H.~Jenkins and J.~Serrin.
\newblock Variational problems of minimal surface type {II}. {B}oundary value
  problems for the minimal surface equation.
\newblock {\em Arch. Rational Mech. Anal.}, 21:321--342, 1966.

\bibitem{ka4}
H.~Karcher.
\newblock Embedded minimal surfaces derived from {S}cherk's examples.
\newblock {\em Manuscripta Math.}, 62:83--114, 1988.


\bibitem{mrt}
L.~Mazet, M. M.~Rodr\'{\i}guez and M.~Traizet.
\newblock Saddle towers with infinitely many ends.
\newblock Preprint.

\bibitem{mazet3}
L.~Mazet.
\newblock A uniqueness result for maximal surfaces in Minkowski $3$-space.
\newblock Preprint.

\bibitem{mr7}
W.~H. Meeks~III and H.~Rosenberg.
\newblock Maximum principles at infinity with applications to minimal and
  constant mean curvature surfaces.
\newblock Preprint.

\bibitem{mrw1}
W.~H. Meeks~III and M.~Wolf.
\newblock Minimal surfaces with the area growth of two planes; the
case of infinite symmetry.
\newblock To appear in J. Amer. Math. Soc.

\bibitem{PeTra1}
J.~P\'{e}rez and M.~Traizet.
\newblock The classification of singly periodic minimal surfaces with genus
  zero and {S}cherk type ends.
\newblock To appear in Trans. Amer. Math. Society.

\bibitem{mrod2}
M.~M. Rodr\'{\i}guez.
\newblock A {J}enkins-{S}errin problem on the strip.
\newblock To appear in J. Geometry and Physics.

\bibitem{mrod1}
M.~M. Rodr\'{\i}guez.
\newblock The space of doubly periodic minimal tori with parallel ends: the
  standard examples.
\newblock To appear in Michigan Math. J.

\bibitem{ros5}
A.~Ros.
\newblock Embedded minimal surfaces: forces, topology and symmetries.
\newblock {\em Calc. Var.}, 4:469--496, 1996.


\bibitem{sche1}
H.~F. Scherk.
\newblock Bemerkungen \"{u}ber die kleinste {F}l\"{a}che innerhalb gegebener
  {G}renzen.
\newblock {\em J. R. Angew. Math.}, 13:185--208, 1835.

\end{thebibliography}

\end{document}